\documentclass[12pt, fleqn]{article}
\usepackage[utf8]{inputenc}
\usepackage{amsmath, amssymb, amsthm}
\usepackage{geometry}
\usepackage{hyperref}
\title{The Relationship Between Euler Numbers and Bernoulli Numbers with Ordered Partitions}
\author{  Kamyar Sepehri Pirayvatloo  \& Kazem Haghnejad Azar
\\University of Mohaghegh Ardabili\\
haghnejad@uma.ac.ir, kamyarpmf1@gmail.com}
\date{}
\begin{document}
\setlength{\mathindent}{1.5em} 
\maketitle

\begin{abstract}
In this paper, for every $n \in \mathbb{N}$, the following relationships between the functions $K_{b}(n)$ and $K_{e}(n)$ and the Bernoulli and Euler numbers are proved:
\begin{align*}
B_{2n} &= - \frac{(2n)!}{2^{2n} - 2} \times K_{b}(n) \\
E_{2n} &= (2n)! \times K_{e}(n)
\end{align*}

where the functions \( K_{b} \) and \( K_{e} \) are defined recursively as follows:

\begin{align*}
K_{b}(0) &= K_{e}(0) = 1 \\
K_{b}(n) &= - \sum_{n' = 0}^{n - 1} \frac{K_{b}( n' )}{( 2( n - n' ) + 1 )!}, \quad n \geq 1 \\
K_{e}(n) &= - \sum_{n' = 0}^{n - 1} \frac{K_{e}( n' )}{( 2(n - n') )!}, \quad n \geq 1
\end{align*}

Furthermore, we present combinatorial interpretations of these functions in terms of ordered partitions of \( n \), as follows:

\begin{align*}
K_{b}(n) &= \sum_{\lambda \vDash n} \frac{( - 1)^{l(\lambda)}}{\prod_{i = 1}^{l(\lambda)} ( 2b_{i} + 1 )!}, \quad n \geq 1 \\
K_{e}(n) &= \sum_{\lambda \vDash n} \frac{( - 1)^{l(\lambda)}}{\prod_{i = 1}^{l(\lambda)} ( 2b_{i} )!}, \quad n \geq 1
\end{align*}

where \( \lambda = (b_{1},b_{2},\ldots,b_{k}) \vDash n \) and \( l(\lambda) = k \).
\end{abstract}

\section{Introduction}
The primary objective of this paper is to develop a unified framework for two fundamental families of numbers, the Bernoulli numbers and the Euler numbers, and to introduce new definitions and computational techniques that arise from this shared structure. By relating these numbers to ordered partitions (compositions) of integers, we reveal a novel connection linking number theory, combinatorics, and special functions. In the sections that follow, we provide the necessary definitions and establish the key relations.

\section{Definitions}
\subsection{Ordered Partitions of n}
An ordered partition (composition) of \( n \) is a sequence of positive integers \( \lambda = (b_{1},b_{2},\ldots,b_{k}) \) such that:

\[
b_{1} + b_{2} + \ldots + b_{k} = n
\]

where the order of the components matters and the number of components \( (k) \) is called the length of the partition.

\subsection{Bernoulli Numbers}
\textbf{Jacob Bernoulli} (1654-1705), a Swiss mathematician and physicist, first systematically studied Bernoulli numbers in 1713. These numbers emerged during his investigations of the sums of powers of natural numbers and became a fundamental tool in number theory and analysis \cite{1}. Over time, the connection of this sequence of numbers with the expansions of many functions, including the zeta function, has been established; thus, Bernoulli numbers hold great importance in mathematics. In general, the sum \( \sum_{k = 1}^{n - 1} k^{r} \) can be written as follows \cite{2}:

\[
\sum_{k = 1}^{n - 1} k^{r} = \sum_{k = 0}^{r} \frac{B_{k} \, r! \, n^{r - k + 1}}{k!(r - k + 1)!}
\]

where \( B_{k} \) is independent of \( r \) and are called Bernoulli numbers.

\subsection{Euler Numbers}
\textbf{Leonhard Euler} (1707-1783) studied Euler numbers in 1755 in his book \textit{Institutiones Calculi Differentialis} \cite{3}. These numbers, which initially appeared in the context of expansions of trigonometric series, gradually found extensive applications in combinatorics and graph theory. By discovering the connection between these numbers and specific permutations, Euler laid foundations that are still used in combinatorial research today. These numbers form a sequence \( E_{n} \) of natural numbers defined by the \textbf{Taylor series} as follows \cite{4}:

\[
\frac{1}{\cosh(t)} = \frac{2}{e^{t} + e^{-t}} = \sum_{n = 0}^{\infty} \frac{E_{n}}{n!} t^{n}
\]

\subsection{Zeta Function}
The Riemann zeta function \( \zeta(s) \) was first studied by \textbf{Leonhard Euler} in the 18th century for real numbers \( (s > 1) \), but \textbf{Bernhard Riemann} (1826-1866) extended it to complex numbers in 1859 \cite{7}.This function is defined for \( \text{Re(s)} > 1 \) as:

\[
\zeta(s) = \sum_{k=1}^{\infty} \frac{1}{k^{s}}
\]

Riemann showed that this function can be \textbf{analytically continued} for all complex numbers except \( s = 1 \) (where it has a simple pole). This function has a deep connection with the distribution of prime numbers \cite{7}.

\subsection{Hurwitz Zeta Function}
The Hurwitz zeta function, a generalization of the Riemann zeta function, was introduced by \textbf{Adolf Hurwitz} (1859-1919) in 1882. For \( \text{Re(s)} > 1 \), it is defined as an absolutely convergent series \cite{8}:

\[
\zeta(s, x) = \sum_{n = 0}^{\infty} \frac{1}{(n + x)^{s}}; \quad \frac{d}{dx} \zeta(s, x) = - s \, \zeta(s + 1, x)
\]

This function can be extended to other values of \( s \) by analytic continuation, except at \( s=1 \) where the function has a simple pole.

Generalized harmonic functions \( H_{x}^{s} \) can be expressed via the Hurwitz zeta function:

\[
H_{x}^{s} = \sum_{n = 1}^{x} \frac{1}{n^{s}} = \sum_{n = 1}^{\infty} \frac{1}{n^{s}} - \sum_{n = 0}^{\infty} \frac{1}{(n + x + 1)^{s}};
\]
\[
H_{x}^{s} = \zeta(s) - \zeta(s, x + 1);
\]
\[
\frac{d}{dx} ( H_{x}^{s} ) = s \, \zeta(s + 1, x + 1)
\]

\subsection{Gamma Function}
The \textbf{Gamma Function} was first introduced by the Swiss mathematician \textbf{Leonhard Euler} in his endeavor to generalize the factorial function to non-integer values. Useful details on the origin and evolution of the Gamma function can be found in \cite{5}. Euler introduced the following analytic function during the years 1729 and 1730:

\[
\Gamma(n) = \int_{0}^{1} (-\ln(t))^{n - 1}  dt, \quad (n>0)
\]

which, by an elementary change of variable, leads to the following equation \cite{5}:

\[
\Gamma(n) = \int_{0}^{\infty} t^{n-1} e^{-t} dt, \quad (n>0)
\]

Connection with the factorial function \cite{5}:

\[
\Gamma(n + 1) = n \Gamma(n) = n!
\]

Some specific values \cite{5}:

\[
\Gamma(1) = 1, \quad \Gamma\left( \frac{1}{2} \right) = \sqrt{\pi}, \quad \Gamma\left( \frac{3}{2} \right) = \frac{\sqrt{\pi}}{2}
\]

Efforts have been made and formulas have been provided that describe the Gamma function. One such formula for \( \text{Re}(x) > - \frac{1}{2} \) in article \cite{6}, page 12, is as follows:

\begin{equation}
\Gamma(x+1) = \frac{ \left( x + \frac{1}{2} \right)^{x + \frac{1}{2}} \sqrt{2\pi} }{ e^{x + \frac{1}{2}} } \times e^{ - \sum_{j=1}^{\infty} \frac{ \zeta(2j, x+1) }{ (2j + 4j^{2}) 2^{2j} } }
\end{equation}

In addition to the Gamma function itself, its logarithmic derivative, the Digamma function, is also very important. To compute this derivative, we proceed as follows.

\[
H_{x}^{2j} = \zeta(2j) - \zeta(2j, x + 1);
\]
\[
\Gamma(x+1) = \frac{ \left( x + \frac{1}{2} \right)^{x + \frac{1}{2}} \sqrt{2\pi} }{ e^{x + \frac{1}{2}} } \times e^{ \sum_{j=1}^{\infty} \frac{ H_{x}^{2j} - \zeta(2j) }{ (2j + 4j^{2}) 2^{2j} } };
\]

From article \cite{9}, page 50, we know:

\begin{equation}
\sum_{j=1}^{\infty} \frac{ \zeta(2j) }{ (2j + 4j^{2}) 2^{2j} } = \ln \sqrt{\frac{\pi}{e}}
\end{equation}

Therefore:

\[
\Gamma(x+1) = \frac{ \left( x + \frac{1}{2} \right)^{x + \frac{1}{2}} \sqrt{2} }{ e^{x} } \times e^{ \sum_{j=1}^{\infty} \frac{ H_{x}^{2j} }{ (2j + 4j^{2}) 2^{2j} } };
\]

We define the function \( A \) to simplify the differentiation of the Gamma function as follows:

\begin{equation}
A = \sum_{j=1}^{\infty} \frac{ H_{x}^{2j} }{ (2j + 4j^{2}) 2^{2j} }
\end{equation}

Therefore:

\begin{equation}
\Gamma(x+1) = \frac{ \left( x + \frac{1}{2} \right)^{x + \frac{1}{2}} \sqrt{2} }{ e^{x - A} };
\end{equation}

\[
\frac{d}{dx} \Gamma(x) = \frac{ \left( x + \frac{1}{2} \right)^{x + \frac{1}{2}} \sqrt{2} }{ e^{x - A} } \left( \ln \left( x + \frac{1}{2} \right) + \frac{dA}{dx} \right)
\]

\subsection{Digamma Function}
The \textbf{Digamma Function} is obtained by differentiating the logarithm of the Gamma function as follows:

\[
\frac{d}{dx} \Gamma(x) = \Gamma(x) \psi(x);
\]
\[
\Gamma(x + 1) \psi(x + 1) = \frac{ \left( x + \frac{1}{2} \right)^{x + \frac{1}{2}} \sqrt{2} }{ e^{x + A(x, - \frac{1}{2})} } \left( \ln \left( x + \frac{1}{2} \right) + \frac{dA}{dx} \right);
\]
\[
\psi(x + 1) = \ln \left( x + \frac{1}{2} \right) + \frac{dA}{dx};
\]
\begin{equation}
\psi(x + 1) = \ln \left( x + \frac{1}{2} \right) + \sum_{n = 1}^{\infty} \frac{ \zeta(2n + 1, x + 1) }{ (2n + 1) 2^{2n} }
\end{equation}

\subsection{Polygamma Functions}
\textbf{Polygamma Functions} are the logarithmic derivatives of the Gamma function, defined as follows:

\[
\frac{d^{m + 1}}{dx^{m + 1}} \ln \Gamma(x) = \psi^{m}(x);
\]
\[
\frac{d \psi(x + 1)}{dx} = \psi^{1}(x + 1);
\]
\[
\psi^{1}(x + 1) = \frac{1}{x + \frac{1}{2}} - \sum_{n = 1}^{\infty} \frac{ \zeta(2n + 2, x + 1) }{ 2^{2n} };
\]
\[
\psi^{2}(x + 1) = \frac{-1}{ \left( x + \frac{1}{2} \right)^{2} } + \sum_{n = 1}^{\infty} \frac{ (2n + 2) \zeta(2n + 3, x + 1) }{ 2^{2n} };
\]
\[
\psi^{3}(x + 1) = \frac{2}{ \left( x + \frac{1}{2} \right)^{3} } - \sum_{n = 1}^{\infty} \frac{ (2n + 2)(2n + 3) \zeta(2n + 4, x + 1) }{ 2^{2n} };
\]

In general, based on differentiation rules, we have:

\begin{equation}
\psi^{y}(x + 1) = ( - 1)^{y - 1} \left( \frac{(y - 1)!}{ \left( x + \frac{1}{2} \right)^{y} } - \sum_{n = 1}^{\infty} \frac{ (2n + y)! \, \zeta(2n + y + 1, x + 1) }{ 2^{2n} (2n + 1)! } \right)
\end{equation}
\begin{equation*}
\psi^{y}(x + 1) = ( - 1)^{y - 1} \left( \frac{(y - 1)!}{ \left( x + \frac{1}{2} \right)^{y} } - \sum_{n = 1}^{\infty} \frac{ (2n + y)! \, \zeta(2n + y + 1, x + 1) }{ 2^{2n} (2n)(2n + 1)(2n - 1)! } \right)
\end{equation*}

Let:

\begin{equation}
f(n) = \frac{1}{2^{2n}(2n)(2n + 1)}
\end{equation}

Then:

\begin{equation}
\psi^{y}(x + 1) = ( - 1)^{y - 1} \left( \frac{(y - 1)!}{ \left( x + \frac{1}{2} \right)^{y} } - \sum_{n = 1}^{\infty} \frac{ f(n) (2n + y)! \, \zeta(2n + y + 1, x + 1) }{ (2n - 1)! } \right)
\end{equation}

\section{Findings and Proofs}

\subsection{Theorem 1}
If we define the functions \( K_{b}(n) \) and \( K_{e}(n) \) as follows for every \( n \in \mathbb{N} \):

\begin{align}
K_{b}(n) &= \sum_{\lambda \vDash n} \frac{( - 1)^{l(\lambda)}}{\prod_{i = 1}^{l(\lambda)} ( 2b_{i} + 1 )!}, \quad n \geq 1 \\
K_{e}(n) &= \sum_{\lambda \vDash n} \frac{( - 1)^{l(\lambda)}}{\prod_{i = 1}^{l(\lambda)} ( 2b_{i} )!}, \quad n \geq 1
\end{align}

where \( \lambda = (b_{1},b_{2},\ldots,b_{k}) \) represents the ordered partitions of \( n \) and \( l(\lambda) = k \) is the length of each partition, then the following relations hold for Bernoulli and Euler numbers:

\begin{align}
B_{2n} &= - \frac{(2n)!}{2^{2n} - 2} \times K_{b}(n) \\
E_{2n} &= (2n)! \times K_{e}(n)
\end{align}

In addition , the functions \( K_{b}(n) \) and \( K_{e}(n) \) can be computed recursively by the following relations:

\begin{equation}
K_{b}(n) = \begin{cases}
1 & n = 0 \\
- \sum_{n' = 0}^{n - 1} \frac{ K_{b}(n') }{ (2(n - n') + 1)! } & n \geq 1
\end{cases}
\end{equation}

\begin{equation}
K_{e}(n) = \begin{cases}
1 & n = 0 \\
- \sum_{n' = 0}^{n - 1} \frac{ K_{e}(n') }{ (2(n - n'))! } & n \geq 1
\end{cases}
\end{equation}

\subsection{Proof}
An equation for the polygamma function exists in article number \cite{11}, page 430, as follows:

\begin{equation}
\psi^{y}(x + 1) = ( - 1)^{y - 1} (y)! \, \zeta(y + 1, x + 1)
\quad \text{where } y \in \mathbb{N},\ x + 1 \notin \mathbb{Z}_0^-
\end{equation}

Therefore, substituting (8) into (15), we have the following.

\[
\zeta(y + 1, x + 1) = \frac{1}{ y \left( x + \frac{1}{2} \right)^{y} } - \sum_{n = 1}^{\infty} \frac{ f(n) (2n + y)! \, \zeta(2n + y + 1, x + 1) }{ (2n - 1)! y! }
\]

By the definition of the Beta function \cite{8}:
\[
B(2n, 2m_{0}) = \frac{(2n-1)!(2m_{0}-1)!}{(2n + 2m_{0} - 1)!}
\]

Now let:
\[
y = 2m_{0} - 1
\]

Then we have the following.

\begin{equation}
\zeta(2m_{0}, x + 1) = \frac{1}{ (2m_{0} - 1) \left( x + \frac{1}{2} \right)^{2m_{0} - 1} } - \sum_{n = 1}^{\infty} \frac{ f(n) \, \zeta( 2(n + m_{0}), x + 1 ) }{ B(2n, 2m_{0}) }
\end{equation}

Then let:

\[
h(m_{0}) = \zeta(2m_{0}, x + 1)
\]
\[
p(m_{0}) = \frac{1}{ (2m_{0} - 1) \left( x + \frac{1}{2} \right)^{2m_{0} - 1} }
\]
\[
M_{k} = m_{0} + \sum_{j = 1}^{k} b_{j}
\]
\[
m_{k} = m_{0} + k
\]
\[
j(a_{k}, b_{k}) = - \frac{ f(b_{k}) }{ B(2b_{k}, 2a_{k}) }
\]
\[
T(M_{k}) = j(m_{0}, b_{1}) j(M_{1}, b_{2}) \times \ldots \times j(M_{k - 2}, b_{k - 1}) \times j(M_{k - 1}, b_{k})
\]

Then:

\[
a_{1} = m_{0} 
\]
\[
a_{k} = M_{k - 1} = m_{0} + \sum_{j = 1}^{k - 1} b_{j} = m_{0} + \sum_{j = 1}^{k - 2} b_{j} +\ B_{k - 1} ;
\]
\[
a_{k} = a_{k - 1} + b_{k - 1}
\]

Therefore, in the function \( T(M_{k}) \), the following two rules hold:

\[
a_{1} = m_{0}, \quad a_{k} = a_{k - 1} + b_{k - 1}
\]

Now we rewrite equation (16) as follows:

\begin{equation}
h(m_{0}) = p(m_{0}) + \sum_{b_{1} = 1}^{\infty} j(m_{0}, b_{1}) h(m_{0} + b_{1});
\end{equation}

\[
h(m_{0}) = p(m_{0}) + \sum_{b_{1} = 1}^{\infty} j(m_{0}, b_{1}) \left( p(m_{0} + b_{1}) + \sum_{b_{2} = 1}^{\infty} j(m_{0} + b_{1}, b_{2}) h(m_{0} + b_{1} + b_{2}) \right);
\]
\[
h(M_{0}) = p(m_{0}) + \sum_{b_{1} = 1}^{\infty} \sum_{b_{2} = 1}^{\infty} j(m_{0}, b_{1}) p(m_{0} + b_{1}) + j(m_{0}, b_{1}) j(m_{0} + b_{1}, b_{2}) h(m_{0} + b_{1} + b_{2});
\]
\[
h(m_{0}) = p(m_{0}) + \sum_{b_{1} = 1}^{\infty} \sum_{b_{2} = 1}^{\infty} j(m_{0}, b_{1}) p(M_{1}) + j(m_{0}, b_{1}) j(M_{1}, b_{2}) h(M_{2});
\]

\begin{equation}
\begin{aligned}
h(m_{0}) = p(m_{0}) &+ \lim_{z \rightarrow \infty} \left( \sum_{b_{1} = 1}^{z} \sum_{b_{2} = 1}^{z} \ldots \sum_{b_{n} = 1}^{z} j(m_{0}, b_{1}) p(M_{1}) \right. \\
&+ j(m_{0}, b_{1}) j(M_{1}, b_{2}) p(M_{2}) + \ldots \\
&+ \left( j(m_{0}, b_{1}) j(M_{1}, b_{2}) \times \ldots \times j(M_{k - 1}, b_{k}) \right) p(M_{k}) + \ldots \\
&\left. + \left( j(m_{0}, b_{1}) j(M_{1}, b_{2}) \times \ldots \times j(M_{z - 1}, b_{z}) \right) p(M_{z}) \right);
\end{aligned}
\end{equation}

\begin{align*}
h(m_{0}) = p(m_{0}) &+ \lim_{z \rightarrow \infty} \left( \sum_{b_{1} = 1}^{z} \sum_{b_{2} = 1}^{z} \ldots \sum_{b_{n} = 1}^{z} T(M_{1}) p(M_{1}) \right. \\
&+ T(M_{2}) p(M_{2}) + \ldots + T(M_{k}) p(M_{k}) + \ldots \\
&\left. + T(M_{z}) p(M_{z}) \right);
\end{align*}

Therefore, the third rule for \( T(M_{k}) \) is established in the above equation as follows:

\[
a_{k} + b_{k} = M_{k}
\]

Hence, in the function \( h(m_{0}) \), the coefficient of \( p(M_{k}) \) is always \( T(M_{k}) \), which follows the aforementioned triple rules, and these rules are the necessary and sufficient condition for the validity of \( T(M_{k}) \) in the equation.

Also, in the above equation, each instance of the function \( p \), such as \( M_{k} \), can take different values like \( m_{z} \), provided that \( z \) is not less than \( k \). More precisely, for every \( \mathcal{A} \geq 1 \), if \( k = z + \mathcal{A} \), then \( p(m_{z}) \neq p(M_{k}) \) because:

\[
M_{k} = m_{0} + \sum_{j = 1}^{k} b_{j} \overset{b \geq 1}{\rightarrow} M_{k} \geq m_{0} + k;
\]
\[
M_{k} \geq M_{0} + z + \mathcal{A};
\]
\[
M_{k} > m_{z}
\]

Therefore, \( M_{z} \) can take the value \( m_{z} \), but \( M_{z + \mathcal{A}} \) will never be equal to \( m_{z} \).

Also, for every \( z > \mathcal{B} \geq 0 \), if \( k = z - \mathcal{B} \), i.e., \( z \geq k \), then:

\[
M_{k} \geq m_{0} + k;
\]
\[
M_{k} \geq m_{0} + z - \mathcal{B};
\]
\[
M_{k} \geq m_{z} - \mathcal{B};
\]

That is, \( M_{z - \mathcal{B}} \), which is a variable, can take the value \( m_{z} \).

From the summary of the two rules above, it can be concluded that the value \( m_{z} \) is only present in \( M_{1} \) to \( M_{z} \), and also \( p(m_{z}) \) is present in \( p(M_{1}) \) to \( p(M_{z}) \), and its coefficients are calculable such that for every \( p(m_{z}) \) there is a coefficient \( T(m_{z}) \) that follows the triple rules. Therefore, the sum of all coefficients of \( p(m_{z}) \) in the function \( h(m_{0}) \) will be the sum of \( T(m_{z}) \)'s that follow the following triple rules:

\[
a_{1} = m_{0}, \quad a_{k} = a_{k - 1} + b_{k - 1}, \quad a_{k} + b_{k} = M_{k}
\]

Suppose we want to calculate the coefficient of \( p(m_{2}) \):

As we showed, \( p(m_{2}) \) must be present in the functions \( p(M_{1}) \) to \( p(M_{2}) \), so we disregard \( p(M_{3}) \) because the minimum value it can have is \( p(m_{3}) \). Therefore:

\[
h(m_{0}) = p(m_{0}) + \sum_{b_{1} = 1}^{\infty} \sum_{b_{2} = 1}^{\infty} j(m_{0}, b_{1}) p(M_{1}) + j(m_{0}, b_{1}) j(M_{1}, b_{2}) p(M_{2});
\]
\begin{align*}
h(m_{0}) = p(m_{0}) &+ j(m_{0}, 1) p(m_{1}) + \sum_{b_{2} = 1}^{\infty} j(m_{0}, 1) j(m_{1}, b_{2}) p(m_{1} + b_{2}) \\
&+ j(m_{0}, 2) p(m_{2}) + \sum_{b_{2} = 1}^{\infty} j(m_{0}, 2) j(m_{2}, b_{2}) p(m_{2} + b_{2});
\end{align*}
which for \( b_2 = 1 \) gives us:
\begin{align*}
h(m_{0}) = p(m_{0}) &+ j(m_{0}, 1) p(m_{1}) + j(m_{0}, 1) j(m_{1}, 1) p(m_{2}) \\
&+ j(m_{0}, 2) p(m_{2}) + j(m_{0}, 2) j(m_{2}, 1) p(m_{3});
\end{align*}
\begin{align*}
h(m_{0}) = p(m_{0}) &+ j(m_{0}, 1) p(m_{1}) + j(m_{0}, 2) j(m_{2}, 1) p(m_{3}) \\
&+ \left( \mathbf{j(m_{0}, 1) j(m_{1}, 1) + j(m_{0}, 2)} \right) p(m_{2})
\end{align*}

\subsection{Result}
1. The function \( h(m_{0}) \) can be computed as an infinite sum of the function \( p \) with the corresponding coefficients.

2. We denote the coefficient of each \( p(m_{z}) \) by \( T(m_{z}) \), which is formed by the product of several functions \( j(a, b) \) such that if we have:

\[
T(m_{z}) = j(a_{1}, b_{1}) \times j(a_{2}, b_{2}) \times \ldots \times j(a_{k - 1}, b_{k - 1}) \times j(a_{k}, b_{k})
\]

then the following triple rules hold:

\[
a_{1} = m_{0}, \quad a_{k} = a_{k - 1} + b_{k - 1}, \quad a_{k} + b_{k} = m_{z} = m_{0} + z
\]

3. If we decompose the function \( h \) according to equation (17) and also decompose the resulting \( h \) functions from the decomposition with the same equation, and factor out the resulting \( p \) functions like \( p(m_{z}) \), and sum their coefficients, i.e., the \( T(m_{z}) \)'s, we will have:

\[
h(m_{0}) = p(m_{0}) + \sum_{z = 1}^{\infty} \left( \sum T(m_{z}) \right) p(m_{z});
\]
\[
h(m_{0}) = p(m_{0}) + \sum_{z = 1}^{\infty} \left( \sum T(m_{0} + z) \right) p(m_{0} + z);
\]

Let:

\[
g(z + m_{0}, m_{0}) = \sum T(m_{0} + z);
\]

Then:

\begin{equation}
h(m_{0}) = p(m_{0}) + \sum_{z = 1}^{\infty} g(z + m_{0}, m_{0}) p(m_{0} + z);
\end{equation}
\begin{equation}
\zeta(2m_{0}, x + 1) = \frac{1}{ (2m_{0} - 1) \left( x + \frac{1}{2} \right)^{2m_{0} - 1} } + \sum_{z = 1}^{\infty} \frac{ g(z + m_{0}, m_{0}) }{ (2m_{0} + 2z - 1) \left( x + \frac{1}{2} \right)^{2m_{0} + 2z - 1} }
\end{equation}

Therefore, we showed that the coefficient of \( p(m_{0} + z) \) in the function \( h(m_{0}) \) is the function \( g(z + m_{0}, m_{0}) \), which is formed by the sum of all \( T(m_{z}) \)'s that follow the triple rules. In the following, we obtain an easier way to compute these coefficients.

\subsection{Lemma 1}
For every \( T(m_{n}) \) from the function \( g(n + m_{0}, m_{0}) \), the following relation holds:

\[
\sum_{j = 1}^{k} b_{j} = n
\]

More precisely, the sum of the \( b \)'s of the functions \( j(a, b) \) in a \( T(m_{n}) \) that satisfies equation (14), i.e., follows the triple rules, is equal to \( n \). That is, if we set \( \lambda = (b_{1}, b_{2}, \ldots, b_{k}) \), then \( \lambda \vDash n \) and \( l(\lambda) = k \).

\subsection{Proof}
\[
T(m_{n}) = j(a_{1}, b_{1}) \times j(a_{2}, b_{2}) \times \ldots \times j(a_{k - 1}, b_{k - 1}) \times j(a_{k}, b_{k})
\]

According to the triple rules:

\[
a_{1} = m_{0}, \quad a_{k - 1} + b_{k - 1} = a_{k}, \quad a_{k} + b_{k} = m_{n} = m_{0} + n
\]

Therefore:

\[
T(m_{n}) = j(a_{1}, b_{1}) \times j(a_{1} + b_{1}, b_{2}) \times \ldots \times j(a_{1} + b_{1} + \ldots + b_{k - 2}, b_{k - 1}) \times j(a_{1} + b_{1} + \ldots + b_{k - 1}, b_{k})
\]
\[
a_{k} = a_{1} + \sum_{j = 1}^{k - 1} b_{j};
\]
\[
a_{k} + b_{k} = a_{1} + \sum_{j = 1}^{k} b_{j} = m_{0} + n;
\]
\[
\sum_{j = 1}^{k} b_{j} = n
\]

Therefore, \( \lambda = (b_{1}, b_{2}, \ldots, b_{k}) \vDash n \) and \( l(\lambda) = k \).

\subsection{Lemma 2}
If we define the function \( K_{b}(n) \) as follows for every \( n \in \mathbb{N} \):

\begin{equation}
K_{b}(n) = \sum_{\lambda \vDash n} \frac{( - 1)^{l(\lambda)}}{\prod_{i = 1}^{l(\lambda)} ( 2b_{i} + 1 )!}, \quad n \geq 1
\end{equation}

Then:

\begin{equation}
g(n + m_{0}, m_{0}) = \frac{(2n - 1)!}{ B(2n, 2m_{0}) \, 2^{2n} } K_{b}(n)
\end{equation}

\subsection{Proof}
According to previous definitions:
\[
f(n) = \frac{1}{2^{2n}(2n)(2n + 1)}
\]
\[
j(a, b) = - \frac{ f(b) }{ B(2b, 2a) } = - \frac{ (2a + 2b - 1)! }{ 2^{2b} (2b + 1)! (2a - 1)! };
\]

\[
T(m_{n}) = j(a_{1}, b_{1}) \times j(a_{1} + b_{1}, b_{2}) \times \ldots \times j(a_{1} + b_{1} + \ldots + a_{k - 2}, b_{k - 1}) \times j(a_{1} + b_{1} + \ldots + a_{k - 1}, b_{k})
\]

Therefore:
\begin{align*}
T(m_{n}) = & - \frac{ (2a_{1} + 2b_{1} - 1)! }{ 2^{2b_{1}} (2b_{1} + 1)! (2a_{1} - 1)! } 
\times - \frac{ (2a_{1} + 2b_{1} + 2b_{2} - 1)! }{ 2^{2b_{2}} (2b_{2} + 1)! (2a_{1} + 2b_{1} - 1)! } \\
& \times \cdots \
 \times - \frac{ (2a_{1} + 2b_{1} + 2b_{2} + \ldots + 2b_{k} - 1)! }{ 2^{2b_{k}} (2b_{k} + 1)! (2a_{1} + 2b_{1} + 2b_{2} + \ldots + 2b_{k - 1} - 1)! };
\end{align*}
\begin{align*}
T(m_{n}) = & - \frac{ (2a_{1} + 2b_{1} - 1)! }{ 2^{2b_{1}} (2b_{1} + 1)! (2a_{1} - 1)! } 
 \times - \frac{ (2a_{1} + 2b_{1} + 2b_{2} - 1)! }{ 2^{2b_{2}} (2b_{2} + 1)! (2a_{1} + 2b_{1} - 1)! } \\
& \times \cdots \
 \times - \frac{ (2a_{1} + 2\sum_{j=1}^{k} b_{j} - 1)! }{ 2^{2b_{k}} (2b_{k} + 1)! (2a_{1} + 2\sum_{j=1}^{k-1} b_{j} - 1)! };
\end{align*}

After canceling common terms in the numerator and denominator, we have:

\[
T(m_{n}) = \frac{ (-1)^{k} (2a_{1} - 1 + 2 \sum_{j = 1}^{k} b_{j})! }{ (2a_{1} - 1)! \, 2^{2 \sum_{j = 1}^{k} b_{j}} \prod_{j = 1}^{k} (2b_{j} + 1)! };
\]
\[
T(m_{n}) = \frac{ (-1)^{k} \times (2a_{1} - 1 + 2n)! }{ (2a_{1} - 1)! \, 2^{2n} \prod_{j = 1}^{k} (2b_{j} + 1)! };
\]
Given that \( a_{1} = m_{0} \), it follows that:
\[
T(m_{n}) = \frac{ (-1)^{k} (2(n + m_{0}) - 1)! }{ (2m_{0} - 1)! \, 2^{2n} \prod_{j = 1}^{k} (2b_{j} + 1)! };
\]
\[
\sum T(m_{n}) = \frac{ (2(n + m_{0}) - 1)! }{ (2m_{0} - 1)! \, 2^{2n} } \times \sum_{\lambda \vDash n} \frac{ (-1)^{l(\lambda)} }{ \prod_{i = 1}^{l(\lambda)} (2b_{i} + 1)! };
\]
\[
\sum T(m_{n}) = \frac{ (2(n + m_{0}) - 1)! }{ (2m_{0} - 1)! \, 2^{2n} } K_{b}(n);
\]
\[
\sum T(m_{n}) = \frac{ (2n - 1)! }{ B(2n, 2m_{0}) \, 2^{2n} } K_{b}(n);
\]
\[
g(n + m_{0}, m_{0}) = \frac{ (2n - 1)! }{ B(2n, 2m_{0}) \, 2^{2n} } K_{b}(n)
\]

Thus, from the ordered partitions of \( n \), we can compute \( K_{b}(n) \) and subsequently \( g(n + m_{0}, m_{0}) \).

\begin{align*}
K_{b}(3) = & \frac{ (-1)^{3} }{ (2(1) + 1)! (2(1) + 1)! (2(1) + 1)! } 
          + \frac{ (-1)^{2} }{ (2(1) + 1)! (2(2) + 1)! } \\
          & + \frac{ (-1)^{2} }{ (2(2) + 1)! (2(1) + 1)! } 
          + \frac{ (-1)^{1} }{ (2(3) + 1)! } \ = - \frac{31}{15120}
\end{align*}
In the following, we will obtain another equation for computing \( K_{b}(n) \).

\subsection{Lemma 3}
For every \( \text{Re}(x) > -\frac{1}{2} \), the following relation holds for the Gamma function:

\begin{equation}
\Gamma\left( x + \frac{1}{2} \right) = \left( \frac{x}{e} \right)^{x} \sqrt{2\pi} \, \times \, \exp\left( - \sum_{n = 1}^{\infty} \frac{a_{n}}{(2n - 1) \left( x + \frac{1}{2} \right)^{(2n - 1)}} \right) 
\end{equation}

where the coefficients \( a_{n} \) are defined as follows:

\begin{equation}
a_{n} = f(n) + \sum_{m_{0} = 1}^{n - 1} g\left( n, m_{0} \right) \, f\left( m_{0} \right) 
\end{equation}

\subsection{Proof}
We begin with the following known relation for the Gamma function from \cite{6}:

\begin{equation*}
\Gamma(x + 1) = \frac{ \left( x + \frac{1}{2} \right)^{x + \frac{1}{2}} \sqrt{2\pi} }{ e^{x + \frac{1}{2}} } \times \exp\left( - \sum_{n = 1}^{\infty} \frac{\zeta(2n, x + 1)}{ (2n + 4n^{2}) \, 2^{2n} } \right) \tag{1}
\end{equation*}

Let:

\begin{equation*}
U = \sum_{m_{0} = 1}^{\infty} \frac{\zeta\left( 2m_{0}, x + 1 \right)}{\left( 2m_{0} + 4m_{0}^{2} \right) \, 2^{2m_{0}}} 
\end{equation*}

\begin{equation*}
h\left( m_{0} \right) = \zeta\left( 2m_{0}, x + 1 \right)
\end{equation*}

\begin{equation*}
f\left( m_{0} \right) = \frac{1}{\left( 2m_{0} + 4m_{0}^{2} \right) \, 2^{2m_{0}}} \tag{7}
\end{equation*}

Thus, we can write:

\begin{equation}
U = \sum_{m_{0} = 1}^{\infty} f\left( m_{0} \right) \, h\left( m_{0} \right) 
\end{equation}

As demonstrated in the previous section of the paper, the function \( h\left( m_{0} \right) \) can be expanded in terms of the functions \( p \) as follows:

\begin{equation*}
h\left( m_{0} \right) = p\left( m_{0} \right) + \sum_{z = 1}^{\infty} g\left( z + m_{0}, m_{0} \right) \, p\left( z + m_{0} \right) \tag{19}
\end{equation*}
\[
h\left( m_{0} \right) = p\left( m_{0} \right) + g\left( 1 + m_{0}, m_{0} \right) \, p\left( 1 + m_{0} \right) + g\left( 2 + m_{0}, m_{0} \right) \, p\left( 2 + m_{0} \right) + \ldots
\]

This expansion possesses the following properties:

\begin{itemize}
    \item For \( n = m_{0} \), the coefficient of \( p(n) \) in the expansion of \( h\left( m_{0} \right) \) is equal to \( 1 \).
    \item For \( n = m_{0} + c \), where \( c \in \mathbb{N} \), the coefficient of \( p(n) \) in the expansion of \( h\left( m_{0} \right) \) is equal to \( g\left( n, m_{0} \right) \).
\end{itemize}

In the following, we will use these properties of the expansion of \( h \) to extract the coefficients of \( p(n) \) from the function \( U \). For this purpose, we first need to determine which \( h \) functions, upon expansion, yield \( p(n) \).

\[
U = \sum_{m_{0} = 1}^{\infty} f\left( m_{0} \right) \, h\left( m_{0} \right) = \sum_{m_{0} = 1}^{n} f\left( m_{0} \right) \, h\left( m_{0} \right) + \sum_{m_{0} = n + 1}^{\infty} f\left( m_{0} \right) \, h\left( m_{0} \right)
\]

According to equation (19), for \( m_{0} > n \), the function \( p(n) \) does not appear in the expansion of \( h\left( m_{0} \right) \), therefore:

\[
U = \sum_{m_{0} = 1}^{n} f\left( m_{0} \right) \, h\left( m_{0} \right) + \text{(sum of other terms without } p(n) \text{)}
\]

Consequently, to find the coefficient of \( p(n) \), we only examine the first part:

\[
U = \sum_{m_{0} = 1}^{n} f\left( m_{0} \right) \, h\left( m_{0} \right) + \ldots 
\]
\[
U =  f(n) \, h(n) + \sum_{m_{0} = 1}^{n - 1} f\left( m_{0} \right) \, h\left( m_{0} \right) + \ldots
\]

Based on the explanations provided earlier:

\begin{enumerate}
    \item The coefficient of \( p(n) \) in the expansion of \( h(n) \) is equal to one. Therefore, one of the coefficients of \( p(n) \) in the function \( U \) will be \( f(n) \).
    \item For every \( n > m_{0} \), the coefficient of \( p(n) \) in the expansion of \( h\left( m_{0} \right) \) is equal to \( g\left( n, m_{0} \right) \). Thus, the other coefficients of \( p(n) \) in the function \( U \), for every \( m_{0} \leq n - 1 \), will be \( f\left( m_{0} \right) \, g\left( n, m_{0} \right) \).
\end{enumerate}

Therefore, the sum of the coefficients of \( p(n) \) in the expansion of \( U \) is equal to:

\begin{equation}
a_{n} = f(n) + \sum_{m_{0} = 1}^{n - 1} g\left( n, m_{0} \right) \, f\left( m_{0} \right) \tag{24}
\end{equation}

Now, having the coefficients \( a_{n} \), we can expand \( U \) in terms of the functions \( p(n) \) as follows:

\begin{equation}
U = \sum_{n = 1}^{\infty} a_{n} \, p(n) 
\end{equation}

Substituting equation (26) into equation (1) yields the following relation for the Gamma function:

\[
\Gamma\left( x + \frac{1}{2} \right) = \left( \frac{x}{e} \right)^{x} \sqrt{2\pi} \, \times \, \exp\left( - \sum_{n = 1}^{\infty} a_{n} \, p(n) \right);
\]
\begin{equation}
\Gamma\left( x + \frac{1}{2} \right) = \left( \frac{x}{e} \right)^{x} \sqrt{2\pi} \, \times \, \exp\left( - \sum_{n = 1}^{\infty} \frac{a_{n}}{(2n - 1) \left( x + \frac{1}{2} \right)^{(2n - 1)}} \right) \tag{23}
\end{equation}

\subsection{Lemma 4}
\begin{equation}
a_{n} = - B(2n, 2m_{0}) g(n + m_{0}, m_{0})
\end{equation}

\subsection{Proof}
We take the derivative from equation (23):

\[
\psi(x + 1) = \ln \left( x + \frac{1}{2} \right) + \sum_{n = 1}^{\infty} \frac{ a_{n} }{ \left( x + \frac{1}{2} \right)^{2n} };
\]
\[
\psi^{y}(x + 1) = (-1)^{y - 1} \left( \frac{ (y - 1)! }{ \left( x + \frac{1}{2} \right)^{y} } - \sum_{n = 1}^{\infty} \frac{ (2n + y - 1)! a_{n} }{ (2n - 1)! \left( x + \frac{1}{2} \right)^{2n + y} } \right);
\]
\[
\zeta(y + 1, x + 1) = \frac{1}{ y \left( x + \frac{1}{2} \right)^{y} } - \sum_{n = 1}^{\infty} \frac{ (2n + y - 1)! a_{n} }{ (2n - 1)! y! \left( x + \frac{1}{2} \right)^{2n + y} };
\]

Let:

\[
y = 2m_{0} - 1
\]

Then:

\begin{equation}
\zeta(2m_{0}, x + 1) = \frac{1}{ (2m_{0} - 1) \left( x + \frac{1}{2} \right)^{2m_{0} - 1} } - \sum_{n = 1}^{\infty} \frac{ \frac{ a_{n} }{ B(2n, 2m_{0}) } }{ (2m_{0} + 2n - 1) \left( x + \frac{1}{2} \right)^{2m_{0} + 2n - 1} }
\end{equation}

Therefore, by comparing the two equations (19) and (28), we obtain:

\begin{equation*}
a_{n} = - B(2n, 2m_{0}) .g(n + m_{0}, m_{0})
\end{equation*}

In the following, we obtain a recursive relation for \( a_{n} \) through this equation.

\subsection{Lemma 5}
\begin{equation}
a_{n} = \frac{1}{ 2^{2n} (2n + 4n^{2}) } - \sum_{n' = 1}^{n - 1} \frac{ \binom{2n - 1}{2n'} }{ 2^{2n'} (2n' + 1) } a_{n - n'}
\end{equation}

\subsection{Proof}
We previously showed:

\begin{equation}
a_{n} = f(n) + \sum_{n' = 1}^{n - 1} g(n, n') f(n')
\end{equation}

\begin{equation}
a_{n} = - B(2n, 2n') g(n + n', n')
\end{equation}

Therefore, the \( n \)-th coefficient can be computed by having the previous coefficients:

\[
a_{n - m} = - B(2n - 2m, 2n') g(n - m + n', n') \overset{m = n'}{\rightarrow}
\]
\[
a_{n - n'} = - B(2n - 2n', 2n') g(n, n');
\]
\[
g(n, n') = - \frac{ a_{n - n'} }{ B(2n - 2n', 2n') } = - a_{n - n'} \frac{ \Gamma(2n) }{ \Gamma(2n - 2n') \Gamma(2n') };
\]
\[
g(n, n') = - a_{n - n'} \frac{ (2n - 1)! }{ (2n - 2n' - 1)! (2n' - 1)! };
\]
\[
g(n, n') = - (2n') \binom{2n - 1}{2n'} a_{n - n'};
\]
\[
a_{n} = f(n) - \sum_{n' = 1}^{n - 1} f(n') (2n') \binom{2n - 1}{2n'} a_{n - n'}
\]
\begin{equation*}
\mathbf{a_{n}} = \frac{\mathbf{1}}{\mathbf{2}^{2n} (2n + 4n^{2})} - \sum_{n' = 1}^{n - 1} \frac{ \binom{2n - 1}{2n'} }{ 2^{2n'} (2n' + 1) } \mathbf{a_{n - n'}} \tag{29}
\end{equation*}

In the following, from the combination of relation number (22) and (27), we obtain:

\begin{equation}
a_{n} = \frac{ - (2n - 1)! }{ 2^{2n} } K_{b}(n);
\end{equation}

Finally, by substituting the above relation into relation (29), a simpler and more optimized method for computing \( K_{b}(n) \) is obtained:

\[
K_{b}(n) = - \frac{1}{(2n + 1)!} - \sum_{n' = 1}^{n - 1} \frac{ K_{b}(n') }{ (2(n - n') + 1)! }, \quad n \geq 1
\]

Now, given that \( K_{b}(1) = - \frac{1}{6} \), we can modify the equation as follows:

\begin{equation*}
K_{b}(n) = \begin{cases}
1 & n = 0 \\
- \sum_{n' = 0}^{n - 1} \frac{ K_{b}(n') }{ (2(n - n') + 1)! } & n \geq 1
\end{cases} \tag{13}
\end{equation*}

\[
\frac{ K_{b}(n) }{ 1! } = - \frac{ K_{b}(n - 1) }{ 3! } - \frac{ K_{b}(n - 2) }{ 5! } - \ldots - \frac{ K_{b}(1) }{ (2n - 1)! } - \frac{ K_{b}(0) }{ (2n + 1)! }
\]
\[
\frac{ K_{b}(0) }{ 1! } = 1
\]
\[
\frac{ K_{b}(1) }{ 1! } = - \frac{ K_{b}(0) }{ 3! }
\]
\[
\frac{ K_{b}(2) }{ 1! } = - \frac{ K_{b}(1) }{ 3! } - \frac{ K_{b}(0) }{ 5! }
\]
\[
\frac{ K_{b}(3) }{ 1! } = - \frac{ K_{b}(2) }{ 3! } - \frac{ K_{b}(1) }{ 5! } - \frac{ K_{b}(0) }{ 7! }
\]

\subsection{Connection with Bernoulli Numbers}
We previously showed:

\[
\psi(x + 1) = \ln \left( x + \frac{1}{2} \right) + \sum_{n = 1}^{\infty} \frac{ a_{n} }{ \left( x + \frac{1}{2} \right)^{2n} };
\]

In article \textbf{\cite{10}}, a similar equation for the digamma function is discussed as follows:

\[
\psi(x + 1) = \ln \left( x + \frac{1}{2} \right) + \sum_{n = 1}^{\infty} \frac{ B_{2n} \left( 1 - \frac{1}{2^{2n - 1}} \right) }{ 2n \left( x + \frac{1}{2} \right)^{2n} };
\]

Therefore, by comparing the two equations above, the relation between the Bernoulli numbers and the coefficients $a_n$ is obtained as follows:

\begin{equation}
a_{n} = \frac{ B_{2n} }{ 2n } \left( 1 - \frac{1}{2^{2n - 1}} \right);
\end{equation}
\[
B_{2n} = \frac{ 2n \times a_{n} }{ \left( 1 - \frac{1}{2^{2n - 1}} \right) }
\]

Also, we previously showed:

\[
a_{n} = \frac{ - (2n - 1)! }{ 2^{2n} } K_{b}(n)\tag{32}
\]

As a result:

\begin{equation}
\mathbf{B_{2n}} = - \frac{ (2n)! }{ 2^{2n} - 2 } \times \mathbf{K_{b}(n)}
\end{equation}

Therefore, using the coefficients \( \mathbf{K_{b}(n)} \), we will be able to easily compute Bernoulli numbers.

\renewcommand{\arraystretch}{1.5}
\begin{center}
\begin{tabular}{|c|c|c|c|c|c|c|}
\hline
\( n \) & 1 & 2 & 3 & 4 & 5 & 6 \\
\hline
\( K_{b}(n) \) & \( -\frac{1}{6} \) & \( \frac{7}{360} \) & \( -\frac{31}{15120} \) & \( \frac{127}{604800} \) & \( -\frac{73}{3421440} \) & \( \frac{1414477}{653837184000} \) \\
\hline
\end{tabular}
\end{center}

\subsection{Connection with Euler Numbers}

First, we define the function \( K_{e}(n) \) as the following determinant:

\begin{equation}
K_{e}(n) =
\begin{vmatrix}
1 & 0 & 0 & \cdots & 1 \\
\frac{1}{2!} & 1 & 0 & \cdots & 0 \\
\frac{1}{4!} & \frac{1}{2!} & 1 & \cdots & 0 \\
\vdots & \vdots & \vdots & \ddots & \vdots \\
\frac{1}{(2n)!} & \frac{1}{(2n - 2)!} & \frac{1}{(2n - 4)!} & \cdots & \frac{1}{2!} & 0
\end{vmatrix}
\tag{34}
\end{equation}

In article \cite{4}, the following relations for computing Bernoulli and Euler numbers based on determinants are provided:

\[
B_{2n} = - \frac{(2n)!}{2^{2n} - 2} \times
\begin{vmatrix}
1 & 0 & 0 & \cdots & 1 \\
\frac{1}{3!} & 1 & 0 & \cdots & 0 \\
\frac{1}{5!} & \frac{1}{3!} & 1 & \cdots & 0 \\
\vdots & \vdots & \vdots & \ddots & \vdots \\
\frac{1}{(2n + 1)!} & \frac{1}{(2n - 1)!} & \frac{1}{(2n - 3)!} & \cdots & \frac{1}{3!} & 0
\end{vmatrix}
\]

\[
E_{2n} = (2n)! \times
\begin{vmatrix}
1 & 0 & 0 & \cdots & 1 \\
\frac{1}{2!} & 1 & 0 & \cdots & 0 \\
\frac{1}{4!} & \frac{1}{2!} & 1 & \cdots & 0 \\
\vdots & \vdots & \vdots & \ddots & \vdots \\
\frac{1}{(2n)!} & \frac{1}{(2n - 2)!} & \frac{1}{(2n - 4)!} & \cdots & \frac{1}{2!} & 0
\end{vmatrix}
\]

Therefore, according to the definition of \( K_{e}(n) \) in equation (34), the Euler relation can be rewritten as follows:

\begin{equation}
E_{2n} = (2n)! \times K_{e}(n) \tag{35}
\end{equation}

Also, substituting equation (33) into the formula for Bernoulli numbers, the determinant representation of \( K_{b}(n) \) is obtained as follows:

\[
K_{b}(n) =
\begin{vmatrix}
1 & 0 & 0 & \cdots & 1 \\
\frac{1}{3!} & 1 & 0 & \cdots & 0 \\
\frac{1}{5!} & \frac{1}{3!} & 1 & \cdots & 0 \\
\vdots & \vdots & \vdots & \ddots & \vdots \\
\frac{1}{(2n + 1)!} & \frac{1}{(2n - 1)!} & \frac{1}{(2n - 3)!} & \cdots & \frac{1}{3!} & 0
\end{vmatrix}
\]

Comparing the two determinants \( K_{b}(n) \) and \( K_{e}(n) \) shows that they share the same structure, the only difference being that the odd factorials (from \( 3! \) to \( (2n+1)! \)) in \( K_{b}(n) \) are replaced by even factorials (from \( 2! \) to \( (2n)! \)) in \( K_{e}(n) \).

Consequently, if the function \( F \) is defined as follows:

\[
F((2i + 1)!) = K_{b}(n) \quad \quad n \geq i
\]

Then:

\[
F\left( (2i)! \right) = K_{e}(n) \quad \quad n \geq i
\]

Therefore, considering relations (9) and (13) previously obtained for \( K_{b}(n) \), and by applying the modifications related to replacing odd factorials with even factorials, the equivalent relations for \( K_{e}(n) \) can be deduced as follows:

\begin{equation}
K_{e}(n) = - \sum_{n' = 0}^{n - 1} \frac{ K_{e}( n' ) }{ ( 2(n - n') )! } \quad \quad n \geq 1 \tag{36}
\end{equation}

\begin{equation}
K_{e}(n) = \sum_{\lambda \vDash n} \frac{( - 1)^{l(\lambda)}}{\prod_{i = 1}^{l(\lambda)} ( 2b_{i} )!} \quad \quad n \geq 1 \tag{37}
\end{equation}

These relations provide independent computational methods for \( K_{e}(n) \) that can be used for the efficient computation of Euler numbers.

\section{Conclusion}

In this paper, we have established a unified framework that connects Bernoulli numbers and Euler numbers to ordered partitions of integers through the functions \( K_b(n) \) and \( K_e(n) \). By deriving recursive relations and combinatorial interpretations for these functions, we have developed novel methods for computing these important sequences, which originate from expansions of special functions such as the zeta, Gamma, and polygamma functions.

This connection not only reveals a deeper interplay between number theory and combinatorics, but also offers practical advantages for numerical evaluations. Representing these numbers through ordered partitions provides a better understanding of their structure and enables more efficient computations. This approach may pave the way for future research in generalizing these relationships to other numerical sequences and their applications in analytic and combinatorial number theory.

\end{document}